\input amsppt.sty
\loadbold
\TagsOnRight
\hsize 30pc
\vsize 47pc
\def\nmb#1#2{#2}         
\def\cit#1#2{\ifx#1!\cite{#2}\else#2\fi} 
\def\totoc{}             
\def\ign#1{}             

\redefine\g{\frak g}

\redefine\a{\frak a}

\redefine\su{\frak su}
\redefine\so{\frak so}

\redefine\k{\frak k}
\redefine\h{\frak h}

\redefine\m{\frak m}
\redefine\n{\frak n}
\redefine\p{\frak p}

\redefine\1{{1\over 2}}

\define\row#1#2#3{#1_{#2},\ldots,#1_{#3}}

\topmatter
\title Positively curved $7$-dimensional manifolds \endtitle
\author F. Podest\`a \& L. Verdiani \endauthor  
\address University of Florence, Ist. Matematico and Dip. di Matematica 
"U.Dini", Italy\endaddress 
\email podesta\@ udini.math.unifi.it;\quad verdiani\@
udini.math.unifi.it \endemail
\thanks \endthanks 
\abstract 
We deal with seven dimensional compact Riemannian manifolds of positive 
curvature which admit a cohomogeneity one action by a compact Lie 
group $G$. We prove that the manifold is diffeomorphic to a sphere if the 
dimension of the semisimple part of $G$ is bigger than six.
\endabstract

\leftheadtext{\smc }
\rightheadtext{\smc }

\endtopmatter
\document
\subhead\totoc\nmb0{1}.  Introduction \endsubhead
\bigskip
Compact Riemannian manifolds of positive curvature have been studied by 
several authors and, in particular, the homogeneous ones have been 
classified many years ago by Wallach ([Wa]) and by Berard-Bergery 
([BB]). Non homogeneous examples are quite rare and have been found in 
dimension $6$ and $7$ by Eschenburg more recently ([Es1], [Es2]). On 
the other hand it is natural to try to construct new examples on 
manifolds which are not homogeneous but admit a large group of 
symmetries; in this context, Searle ([Se]) has classified 
cohomogeneity one, positively curved Riemannian manifolds in 
dimension $5$ and $6$, finding, up to diffeomorphisms, only spheres 
and complex projective space, while the $4$-dimensional case was treated in [HK]. 
\par
The aim of this note is to study the seven dimensional case, under the 
assumption of a cohomogeneity one action of a compact Lie group of 
isometries. More precisely, we say that a compact Lie group  $G$ acts 
on a compact manifold $M$ by cohomogeneity one if it has a hypersurface 
orbit; in this case, the orbit space $M/G$ is homeomorphic to $S^{1}$
or to a closed interval $[0,1]$. For a detailed exposition and the 
general results, we refer to [Br],[AA],[AA1]. 
\par
Our main result is the following:
\proclaim{Theorem}
Let $M^{7}$ be a compact, positively curved seven dimensional 
Riemannian manifold. Let $G$ be a compact Lie group $G$ 
acting isometrically and almost faithfully on $M^{7}$ by cohomogeneity one. 
If the semisimple part of $G$ has dimension bigger than $6$, then 
$M^{7}$ is diffeomorphic to a sphere $S^{7}$.
\endproclaim
\remark{Remark}
Actually, we prove that, if we want to find a positively curved, 
cohomogeneity one manifold not diffeomorphic to $S^7$, then 
the only "candidate" group is $G= SU(2)\times SU(2)$.
This last case is much more complicated to be handled 
 with and will be object of a forthcoming paper.\par
It is worthy noting here that, for istance, the Berger
space $V^7={{Sp(2)}\over{SU(2)}}$ (see [Be]), endowed with 
the normal homogeneous metric of positive curvature, admits
an isometric $G$-action of cohomogeneity one, where 
$G=Sp(1)\times Sp(1)\subset Sp(2)$. We recall that the space 
$V^7$ is obtained by considering the immersion of $SU(2)$ into 
$Sp(2)$ corresponding to the (essentialy unique) irreducible 
representation $\rho$ of $SU(2)$ on ${\Bbb C}^4\cong {\Bbb H}^2$, which 
is of quaternionic type; the space $V$ is not homeomorphic to $S^7$, 
although it has the same real cohomology (see [Be]).\par
In order to prove that the action of $G$ on $V^7$ has cohomogeneity one, 
we observe the following. If $X\in\a_1$, the Lie algebra of 
$SU(2)$, then $\rho(X)$ leaves exactly two quaternionc lines invariant,
which we call $Q_1(X),Q_2(X)$ and the two restricted representations 
$\rho(X)|_{Q_i(X)}$ are not equivalent. It then follows that, for 
generic $g\in Sp(2)$, the intersection $G\cap gSU(2)g^{-1}$ will be 
finite and a principal $G$-orbit has dimension $6$.\endremark
\bigskip
\subhead\totoc\nmb0{2}. Proof of the main Theorem  \endsubhead
\bigskip
First of all we recall that the manifold $M^{7}$, being of positive 
curvature and compact, has finite fundamental group, so that we shall 
always suppose that $M^{7}$ is simply connected. This implies that 
there is no fibration of $M^{7}$ on $S^{1}$, so that there are exactly 
two singular orbits. Throughout the following, we will use the symbols 
$K$ to denote a fixed regular isotropy subgroup and $H, H'$ to denote 
two singular isotropy subgroups with $K\subset H\cap H'$; it follows 
from the general theory that $H/K$ and $H'/K$ are diffeomorphic to 
spheres of positive dimensions, since there cannot be execeptional 
orbits (see [Br], [AA]).
\par 
We shall use the following results ([GS]):
\bigskip
\proclaim {Theorem 1.1, [GS]} Let ($M^n,g$) be a simply connected,
$n$-dimensional compact CPCM.\roster
\item If a torus $T^k$ acts
isometrically on $M$, then  $k\leq [{{n+1}\over 2}]$, with equality if
and only if $M^n$ is diffeomorphic to a sphere $S^n$ or a complex
projective space ${\Bbb {CP}}^{n/2}$.
\item If $T^1$ or $SU(2)$ acts isometrically on $M^n$ with a fixed point set of 
codimension two or less than four respectively, then $M^n$ is diffeomorphic to 
a sphere $S^n$, to a complex projective space ${\Bbb {CP}}^{n/2}$ or to a 
quaternionic projective space ${\Bbb {HP}}^{n/4}$.\endroster
\endproclaim
\bigskip
So, we will suppose that the group $G$ has rank less or equal to $3$. Furthermore, 
since the regular orbit has dimension $6$, the group $G$ must have rank at least $2$.
\par
We will divide our analysis according to the difference $d$ between 
the rank of $G$ and the rank of the regular isotropy subgroup $K$. Moreover we note 
that, since the group $G$ acts almost effectively 
on $M$, it acts almost effectively on the regular orbit too, so that $K$ cannot contain any ideal
of $G$.\par \bigskip
\proclaim {Case $d=0$}\endproclaim
\medskip
In this case, $K$ has maximal rank and therefore $G$ is semisimple. 
We will now subdivide our study according to the rank of $G$: \par
\proclaim {Subcase rank($G$) = $2$}\endproclaim 
If $G$ is not simple, 
then $G$ is locally isomorphic to $SU(2)\times SU(2)$ and $K$ must be $2$-dimensional, 
so that $G/K$ has dimension $4$. So, we are left with the case $G$ simple. 
Now, all compact rank $2$, simple Lie groups are locally isomorphic to $SU(3),Spin(5),G_2$ (recall
that $Sp(2)=Spin(5)$). We examine each case separately:\par
\noindent a)\ If $G=G_2$, then $K$ has dimension $8$ and has maximal rank. Since the maximal
subalgebras  of maximal rank of $\g_2$ are $\a_2$ and $2\a_1$ (see [GG]), we see that the Lie algebra
$\k$  of $K$ must be maximal, isomorphic to $\a_2$. In this case, the Lie algebra of a singular
isotropy subgroup must coincide with $\g$, since there are no exceptional orbits. It then follows
that the action of $G$ has exactly two fixed points and the manifold is diffeomorphic to $S^7$.\par
\noindent b)\ If $G=Spin(5)$, then $K$ has dimension $4$ and rank $2$; so the only 
possibility for $\k$ is $\k\cong \Bbb R+\a_1$. Again the Lie algebra $\k$ is maximal and $G$ 
should have a fixed point by the same argument as above. But $Spin(5)$ does not 
act transitively on a $6$-dimensional sphere (see [AA]) and this case is ruled out.\par
\noindent c)\ If $G=SU(3)$, then $K$ has dimension $2$ and rank $2$, so that $K^o$, the 
connected component of $K$, coincides with a maximal torus $T^2$ of $G$. We fix $K$ once for all.
We now consider  a singular isotropy subgroup $H$, with Lie algebra $\h$: now 
$H$ contains $K$ and $H/K$ is  diffeomrophic to a sphere. It is not difficult to see 
that the only possibility for $\h$ is  $\h\cong \Bbb R+\a_1$, maximal Lie subalgebra of 
maximal rank in $\a_2=\g$. We observe that,  each singular orbit is of codimension $4$ in 
$M^7$, so that it is simply connected; therefore $H$ is connected, isomorphic to $U(2)$ and $K$ is
also connected, isomorphic to $T^2$. There are
exactly three immersions of $U(2)$ into $SU(3)$, containing the maximal torus $T^2$ and they 
are mutually conjugate by the Weyl group. It then follows that we have to consider exactly two
cases: the first when the two singular isotropy subgroups $H,H'$ are isomorphic to $U(2)$ but with
different immersions and the second one, when $H=H'$.\par
\medskip
\proclaim {Lemma 3.1} Given the triple ($H,K,H'$) of subgroups of $G=SU(3)$ with $K=T^2$, maximal
torus and $H,H'\cong U(2)$, then \roster
\item if $H\not= H'$, then the manifold is diffeomorphic to $S^7$;
\item if $H=H'$, then the manifold does not carry any positively curved, $G$-invariant metric.
\endroster
\endproclaim
\medskip
\demo {Proof} The first case is easily handled, since the $7$-dimensional sphere admits a 
cohomogeneity one action of the group $SU(3)$, induced by the adjoint representation, which 
admits a triple of subgroups as in (1). \par
In case $H=H'$, we decompose $\g=\k+\m_1+\m_2+\m_3$, where $\m_i,\ i=1,2,3$ are two-dimensional,
irreducible and mutually inequivalent $\k$-modules. We can fix $\h=\h'=\k+\m_1$. 
We now fix a non zero vector $v\in
\m_2$ and  consider a normal geodesic $\gamma:\Bbb R\to M$ w.r.t. a positively
curved  $G$-invariant metric $g$ on $M$; we choose $\gamma$ so that it induces the 
triple $\theta=(H,K,H')$. \par
First of all, we claim that the Killing vector field $X$ induced by $v$ on $M$ 
never vanishes along $\gamma$. This is clear since $\h=\h'$ and $v\not\in\h$.\par
We now consider the smooth function $f(t)=||X||_{\gamma(t)}$ for $t\in\Bbb R$
and we claim that $f$ is a concave positive function, which is not possible.
This will conclude our proof. \par
It will be enough to check that $f"(t)<0$ for all $t$ such that $\gamma(t)$ is 
a regular point. First, we observe that the 
tangent space to a regular orbit splits into $K$-irreducible, mutually inequivalent submodules,
so that the shape operator of the regular orbit hypersurface will 
preserve each submodule and will be a multiple of the identity operator on
$\m_2$. Therefore, if we denote by $D$ the Levi-Civita connection of $g$, we 
have that $D_{\gamma(t)'}X$ is a multiple of $X_{\gamma(t)}$; we then have 
$$\eqalign {2R_{X\gamma'X\gamma'} &= 2||D_{\gamma'}X||^2 - {{d^2}\over{dt^2}}
f^2\cr
{}&= 2{{g(D_{\gamma'}X,X)^2}\over{f^2}} - 2(f')^2 - 2ff" = -2ff" > 0,\cr},$$
since $g(D_{\gamma'}X,X) = ff'$. 
\qed\enddemo
\bigskip
We now proceed considering the next\par
\proclaim {Subcase rank($G$) = $3$}\endproclaim 
We have to distinguish $G$ simple or not simple.\par
Let us suppose that $G$ is not simple: then $G$ is locally 
isomorphic to either $G\cong SU(2)^3$ or $G=SU(2)\times G_1$ with $G_1$ simple of rank $2$.\par
If $G\cong SU(2)^3$, then $K^o$ must be of the form $K^o=(T^1)^3$, where each $T^1$ is a maximal
torus in $SU(2)$. But then a singular isotropy subgroup $H$ should have a connected 
component equal to $SU(2)\times (T^1)^2$; both singular orbits $G/H$ and $G/H'$ would then be in
codimension $3$, hence they would be simply connected. Moreover, since the $G$-action on each
singular orbit is not faithful, they would be totally geodesic (see [PV]) and diffeomorphic to 
$S^2\times S^2$; this is not possible, because $S^2\times S^2$ does not carry any
homogeneous positively curved metric (see [HK]).\par
If $G=SU(2)\times G_1$, where $G_1$ is simple of rank $2$, then $K^o$ is of the form $K^o=T^1\times
K_1$, where $K_1\subset G_1$ subgroups of maximal rank; since $\dim G/K = 6$, we have that
$\dim G_1/K_1 = 4$ and hence $\dim G_1\leq 10$. Now $G_1$ can be either $SU(3)$ or $Spin(5)$. \par
If $G_1=SU(3)$, then $G=SU(2)\times SU(3)$ with $K^o=T^1\times U(2)$ and since $U(2)$ is maximal 
in $G_1$, any singular isotopy subgroup $H\supset K$ must be of the form $H^o=SU(2)\times U(2)$;
so any singular orbit is totally geodesic (see [PV]) and of codimension $3$, contradicting Frankel
Theorem (see [Fr]). \par
If $G_1=Spin(5)$, then $\dim \k_1=6$: looking at the list of all maximal subalgebras of maximal rank
in $\so(5)$, we see that $\k_1$ is maximal and siomorphic to $2\a_1$. Again, the same argument as
above rules this case out.\par
We are left with the case where $G$ is simple of rank $3$. Then $G$ is locally isomorphic 
to $SU(4),Spin(7),Sp(3)$. We have that $\dim K=\dim G-6$ and $Sp(3)$ does not have any such
subgroup. \par
In case $G=SU(4)$, we have that $\dim K=9$. But a $9$-dimensional, rank $3$ subalgebra of 
$\a_3$ must be maximal, isomorphic to $\a_2+\Bbb R$; this means that any singular isotropy 
subalgebra must coincide with $\g$. But in this case $G/K$ is not diffeomorphic to a sphere, so 
it is impossible.\par
In case $G=Spin(7)$, the same kind of arguments show that $K^o=Spin(6)$, maximal subalgebra;
then $G$ must have exactly two fixed points and $M$ is diffeomorphic to $S^7$.
\bigskip
\proclaim {Case $d=1$}\endproclaim
\bigskip
We subdivide this case into two subcases, according to the rank of $G$ equal to $2$ o $3$.\par
\noindent (a)\ If rank($G$)=$2$, then the rank of $K$ is $1$, hence $K^o\cong T^1$ or $SU(2)$. 
Therefore $\dim G=7$ or $9$. But there is no compact group $G$ of rank $2$ and dimension 
$7$ or $9$.\par
\noindent (b)\ If the rank of $G$ is $3$, then $K^o$ belongs to the 
list:
$$\{T^2,T^1\times SU(2), SU(2)^2, SU(3), Spin(5), G_2\}\ .$$
It then follows that $\dim K$ belongs to the set $\{2,4,6,8,10,14\}$ and therefore
$\dim G\in \{8,10,12,14,16,20\}$. Now it easy to see that there is no compact group of
rank $3$ and with the indicated dimension.
\bigskip
\proclaim {Case $d=2$}\endproclaim
\bigskip
If the rank of $G$ is $2$, then $K$ is discrete and $G$ has dimension $6$, hence $G$ 
is locally isomorphic to $SU(2)^2$.\par
If the rank of $G$ is $3$, then $K^o$ is either $T^1$ or $SU(2)$. Therefore $\dim G\in\{7,9\}$.
So $G$ cannot be simple. If $G$ is semisimple, not simple, then the only possibility is $G\cong
SU(2)^3$, while, if $G$ is not semisimple, then we have two possibilities, namely 
$G\cong T^1\times SU(2)^2$ or $T^1\times SU(3)$.\par
\bigskip
Summing up, we found that in cases $d=0$ and $d=1$, the manifold $M$ must be diffeomorphic to 
the sphere $S^7$. If $d=2$, then we have the following possibilities for the pair ($G,K^o$):
$$\vbox{\offinterlineskip  \halign {\strut\vrule \hfil \ $#$\ \hfil
&\vrule\hfil \  $#$\ \hfil  &\vrule \hfil\ $#$\ \hfil 
\vrule\cr 
\noalign{\hrule } 
\ {n.}_{\phantom {\sum_1^N}} &  
 \ {G}_{\phantom {\sum_1^N}} &  
 \ {K^o}_{\phantom {\sum_1^N}} 
 \cr 
\noalign{\hrule depth 1 pt} 
\ {1} &
 {T^1\times SU(3)}
 & \ {SU(2)} \cr 
\noalign{\hrule}
\ {2} &
 {SU(2)^3}
 & \ {SU(2)} \cr 
\noalign{\hrule}
\ {3} &
 {T^1\times SU(2)^2}
 & \ {T^1}\cr 
\noalign{\hrule}
\ {4} &
 {SU(2)^2}
 & \ {\{1\}} \cr 
\noalign{\hrule}}}$$
\centerline {{\text {Table 1.}}}
We now prove the following 
\proclaim {Lemma 2.2} Case (1) in Table 1 occurs only if $M$ is diffeomorphic to
$S^7$.\endproclaim
\bigskip
\demo {Proof} First of all we have to identify the subgroup $K^o\cong SU(2)$ inside 
$G$. Note that there are, up to conjugation, exactly two immersions of $SU(2)$ into 
$SU(3)$, corresponding to an irreducible or reducible representation of $SU(2)$ on 
$\Bbb C^3$. We want to prove that $SU(2)$ inside $SU(3)$ must correspond to a reducible 
representation. Indeed, if not, we consider the Lie algebra $\h$ of a singular isotropy 
subgroup $H\supset K$: since $H/K$ must be diffeomorphic to a sphere and $\su(2)$ acts irreducibly
on $\Bbb C^3$, then the only possibility for $\h$ is $\h\cong \su(2)+\Bbb R$. Moreover, since
$\k$ has trivial centralizer in $\su(3)$, we get that $\h$ contains the center of $\g$.
So, the action of $G$ on $G/H$ is not faithful and any singular orbit is totally geodesic
in codimension $2$; but this contradicts Frankel Theorem (see [Fr]).\par
We then have that the centralizer of $\k$ in $\g$ has real dimension $2$ and therefore 
$K^o=SU(2)$ fixes $3$ dimensions; but then by Theorem 1.1, the manifold is diffeomorphic 
to $S^7$. \qed\enddemo
\bigskip
\proclaim {Lemma 2.3} Case (2) in Table 1 occurs only if $M$ is diffeomorphic to $S^7$.
\endproclaim
\bigskip
\demo {Proof} Since the Lie algebra $\k$ cannot contain ideals of $\g$, we have only two
possibilities, namely ($i$): $\k=\su(2)^\Delta\subset 2\su(2)$ or 
($ii$): $\k=\su(2)^\Delta\subset 3\su(2)$,
where the symbol $\Delta$ means diagonal embedding.\par
In case ($i$), $\k$ fixes exactly four dimension and Theorem 1.1 implies the claim.\par
In case ($ii$), we consider a singular isotropy subalgebra $\h$; since $\k$ has trivial centralizer
in $\g$, we see that $\h\cong \su(2)+\su(2)$ containing $\k$ diagonally. But then, $\h$
must contain an ideal $\su(2)$ of $\g$; if we call $\p$ such an ideal, we see that 
the Lie group $P=SU(2)$ corresponding to $\p$ fixes pointwise the orbit $G/H$ in codimension 
four and again by Thm. 1.1 we get our claim.\qed\enddemo
\bigskip
We now want to analyze the case (3) carefully. In this case $G=T^1\times SU(2)^2$ and 
the connected component $K^o=T^1$. We look for possible singular isotropy subalgebras
$\h\supset\k$: if $\n$ denotes the kernel of the slice representation of $\h$, then 
we have the following possibilities: if $\n=\k$, then $\h$ can be 
isomorphic to 
$(i): \Bbb R+\su(2)$ or 
$(ii): \Bbb R^2$ and if $\n$ is trivial, then $\h$ can be isomorphic 
to 
$(iii): \Bbb R+\su(2)$ or 
$(iv): \su(2)$.\par
We now observe the following:\roster
\item  the Lie algebra $\k$ can be supposed to have non trivial projection onto each
$\su(2)$-factors; otherwise it is easy to show that $\k$ would fix $5$ dimensions and 
we could apply Thm.1.1. It then follows that the centralizer of $\k$ in $\g$ is abelian 
and this excludes the possibility (i) for $\h$.
\item If $\h$ is isomorphic to $\Bbb R+\su(2)$, then the semisimple part of $\h$ can be supposed to
be immersed diagonally, otherwise it would fix three dimensions and again we could apply Thm 1.1.
So the center of $\h$ would coincide with the center of $\g$; then we could restrict the action 
of $G$ to its semisimple part $G^s=SU(2)^2$, which would still act transitively on the singular
orbit $G/H$ and the semisimple part $\h^s$ would act by cohomogeneity one on the normal space 
to $G/H$. Therefore we would reduce to case (1) in Table 1. 
\item If $\h$ is isomorphic to $\su(2)$, then again we can suppose it to be diagonally 
embedded. The tangent space of $G/H$ at some point $p$ is an $\su(2)$-module and its 
complexification splits, as $\h$-module, as $\Bbb C + S^2(\Bbb C^2)$. Now, the second 
fundamental form $h$ of $G/H$ at $p$ gives rise to an $\su(2)$-fixed vector in the space 
$S^2(\Bbb C + S^2(\Bbb C^2))\otimes S^2(\Bbb C^2)$, which is easily seen to have no such vector. 
Therefore, $G/H$ is totally geodesic and finitely covered by $S^1\times S^3$: this is impossible,
since, being positively curved, it should have finite $\pi_1$.\endroster
It the follows that both singular isotropy subalgebras can be supposed to be isomorphic to 
$\Bbb R^2$. It is also clear that if $\h$ or $\h'$ contains the center of $\g$, then we could 
apply Thm 1.1; therefore we can suppose that neither $\h$ nor $\h'$ contains the center. \par
We now consider the decomposition $\g=\h+\Bbb R+\m_1+\m_2$, where $\Bbb R$ is a trivial $\h$-module
and $\m_i,\ i=1,2$ are irreducible $\h$-modules (note that we are supposing that $\h$ does not
contain the center and that $\k$ is embedded diagonally) of real dimension two. \par
It is easy to see that $\m_1$ and $\m_2$ are $\h$-inequivalent modules, otherwise $\h$ would
contain the center of $\g$. We now consider the second fundamental form $h$ of $G/H$: since 
$\h$ acts trivially on $\Bbb R$, it follows that $\Bbb R\subset \ker h$. 
Moreover we consider
the kernels $\n_i,\ i=1,2$ of the actions of $\h$ on $\m_i$: it is clear that each $\n_i$ is not 
equal to $\k$ and therefore it acts transitively on the normal space. It follows that 
$h|_{\m_i\times \m_i}=0$. Moreover $\n_1$ does not act trivially on $\m_2$ 
(otherwise $\n_1$ should
coincide with the center of $\g$, which is not contained in $\g$): therefore, by the invariance 
of $h$ under the $\h$-action, we see that, if $h|_{\m_1\times\m_2}\not= 0$, 
then $\m_2$ is equivalent
to $V^*$, where $V$ denotes the normal space to $G/H$. The same argument with $\n_2$, shows then that
$\m_1$ and $\m_2$ are equivalent, a contradiction. Therefore $h=0$. 
But then $G/H$ is a positively curved, homogeneous manifold of dimension $5$, hence finitely covered 
by $S^5$; on the other hand the group $G$ cannot act transitively on a $5$-dimensional sphere (see
e.g. [AA]).\par
We have therefore proved that case (3) can be reduced to case (4) in Table 1, if we want to 
discard manifolds which are diffeomorphic to spheres. This concludes the proof of our main theorem.
\qed

\enddocument

\bigskip\bigskip\bigskip  
\Refs  
\widestnumber\key{AAAA}  
\ref  
\key AA  
\by A.V. Alekseevsky and D.V.Alekseevsky  
\paper G-manifolds with one dimensional orbit space  
\jour Adv. in Sov. Math.  
\vol 8  
\yr 1992   
\pages 1--31  
\endref  
\ref  
\key AA1  
\bysame   
\paper Riemannian G-manifolds with one dimensional orbit space  
\jour Ann. Glob. Anal. and Geom.  
\vol 11  
\yr 1993  
\pages 197--211  
\endref 
\ref 
\key BB
\by L. Berard-Bergery
\paper Les vari\'et\'es Riemanniennes homog\`enes simplement connexes 
de dimension impaire \`a courbure strictement positive
\jour J. Math. Pures Appl.
\vol 55
\yr 1976
\pages 47--68
\endref
\ref
\key Be
\by M. Berger
\paper Les vari\'et\'es Riemanniennes homog\`enes normal simplement connexes 
\`a courbure strictement positive 
\jour Ann. Scuola Normale Superiore, Pisa
\vol 15
\yr 1961
\pages 179--246
\endref
\ref  
\key Br  
\by G.E. Bredon  
\book Introduction to compact transformation groups  
\publ Acad. Press N.Y. London  
\yr 1972  
\endref  
\ref
\key Es1
\by J.H. Eschenburg
\paper New examples of manifolds with strictly positive curvature
\jour Invent. Math.
\vol 66
\yr 1982
\pages 469--480
\endref
\ref
\key Es2
\bysame
\paper Inhomogeneous spaces of positive curvature
\jour Diff. Geom. and Its Appl. 
\vol 2
\yr 1992
\pages 123--132
\endref
\ref
\key Fr
\by T.T Frankel
\paper Manifolds of positive curvature
\jour Pacific J. Math. 
\vol 61
\yr 1961
\pages 165--174
\endref
\ref
\key GG
\by M. Goto, F.D. Grosshans
\book Semisimple Lie Algebras
\publ Lecture Notes in Pure and Appl. Math., Marcel Dekker
\vol 38
\yr 1978
\endref
\ref
\key GS
\by K. Grove and C. Searle
\paper Positively curved manifolds with maximal symmetry-rank 
\jour J. Pure Appl. Algebra
\vol 91
\yr 1994
\pages 137--142
\endref
\ref
\key HK
\by W.Y. Hsiang and B. Kleiner 
\paper On the topology of positively curved 4-manifolds with symmetry
\jour J. Diff. Geometry 
\vol 30
\yr 1989
\pages 615--621
\endref
\ref
\key PV
\by F. Podest\`a, L. Verdiani
\paper Totally geodesic orbits of isometries 
\publ preprint
\yr 1997
\endref
\ref
\key Se
\by C. Searle
\paper Cohomogeneity and positive curvature in low dimension
\jour Math.Z.
\vol 214
\yr 1993
\pages 491--498
\endref
\ref
\key Wa
\by N.R. Wallach 
\paper Compact Homogeneous manifolds with strictly positive curvature
\jour Ann. of Math. 
\vol 96
\yr 1972
\pages 277--295
\endref

\endRefs
\bye